\documentclass{article}

\usepackage{longtable}
\usepackage{booktabs}
\usepackage{charter}
\usepackage{geometry}
\usepackage{graphicx}
\usepackage{amsmath,amssymb}
\usepackage{xcolor}
\usepackage{youngtab}
\usepackage{pifont}
\usepackage{enumitem}

\usepackage[colorlinks=true,bookmarks=false]{hyperref}
\definecolor{myurlcolor}{rgb}{0.6,0,0}
\definecolor{mycitecolor}{rgb}{0,0,0.8}
\definecolor{myrefcolor}{rgb}{0,0,0.8}
\hypersetup{linkcolor=myrefcolor,citecolor=mycitecolor,urlcolor=myurlcolor}

\definecolor{lblue}{HTML}{cce6ff}
\definecolor{lyellow}{HTML}{e0ba4f}
\definecolor{tan}{HTML}{fcba03}
\definecolor{purple}{HTML}{bd358f}
\definecolor{lightgreen}{HTML}{7fde31}
\definecolor{darkorange}{HTML}{b55a2a}
\definecolor{darkgreen}{HTML}{4e8f18}
\definecolor{lightblue}{rgb}{0.9,0.95,1}
\definecolor{darkblue}{rgb}{0,0,100}
\definecolor{brickred}{rgb}{0.7,0,0}

\newcommand{\define}[1]{{\bf \boldmath{#1}}}

\newcommand\Z{{\mathbb Z}}   
\newcommand\R{{\mathbb R}}
\newcommand\C{{\mathbb C}}

\renewcommand\H{{\mathbb H}}
\newcommand{\Oct}{\mathbb{O}}

\newcommand{\U}{{\mathrm U }}
\newcommand{\SU}{{\mathrm{SU}}}
\newcommand{\GL}{{\mathrm{GL}}}
\newcommand{\SL}{{\mathrm{SL}}}
\renewcommand{\O}{{\mathrm{O}}}
\newcommand{\SO}{{\mathrm{SO}}}
\newcommand{\Sp}{{\mathrm{Sp}}}

\newcommand{\g}{\mathfrak{g}}
\newcommand{\ad}{\mathrm{ad}}
\newcommand{\tr}{\mathrm{tr}}

\newcommand{\maps}{\colon}

\newcommand\sbullet{\mathbin{\hbox{\raisebox{-0.75em}{\scalebox{1.5}{$\bullet$}}}}}

\renewcommand{\texttt}[1]{%
  \begingroup
  \ttfamily
  \begingroup\lccode`~=`/\lowercase{\endgroup\def~}{/\discretionary{}{}{}}%
  \begingroup\lccode`~=`[\lowercase{\endgroup\def~}{[\discretionary{}{}{}}%
  \begingroup\lccode`~=`.\lowercase{\endgroup\def~}{.\discretionary{}{}{}}%
  \catcode`/=\active\catcode`[=\active\catcode`.=\active
  \scantokens{#1\noexpand}%
  \endgroup
}

\usepackage{tikz}
\usepackage{tikz-cd}
\usetikzlibrary{matrix,arrows,decorations.pathmorphing,positioning}
\usetikzlibrary{intersections,decorations.markings}
\usetikzlibrary{arrows,positioning,fit,matrix,shapes.geometric,external}
\usetikzlibrary{backgrounds,circuits,circuits.ee.IEC,shapes,fit,matrix}

\tikzstyle{inarrow}=[->, >=stealth, shorten >=.03cm,line width=1.5]

\makeatletter
\newcommand{\xRightarrow}[2][]{\ext@arrow 0359\Rightarrowfill@{#1}{#2}}
\makeatother

\begin{document}

\thispagestyle{empty}

  \begin{center}
    {\Huge\textbf{Coxeter and Dynkin Diagrams}}  
 \\[4em]
  {John Baez}
  \\
  {Department of Mathematics}
  \\
  {University of California, Riverside}
  \\[0.5em]{{based on ``week62--week65'' and ``week182'' of \\ \emph{\href{http://math.ucr.edu/home/baez/TWF.html}{This Week's Finds}}}}
  \\[3em]
  \includegraphics[scale = 0.17]{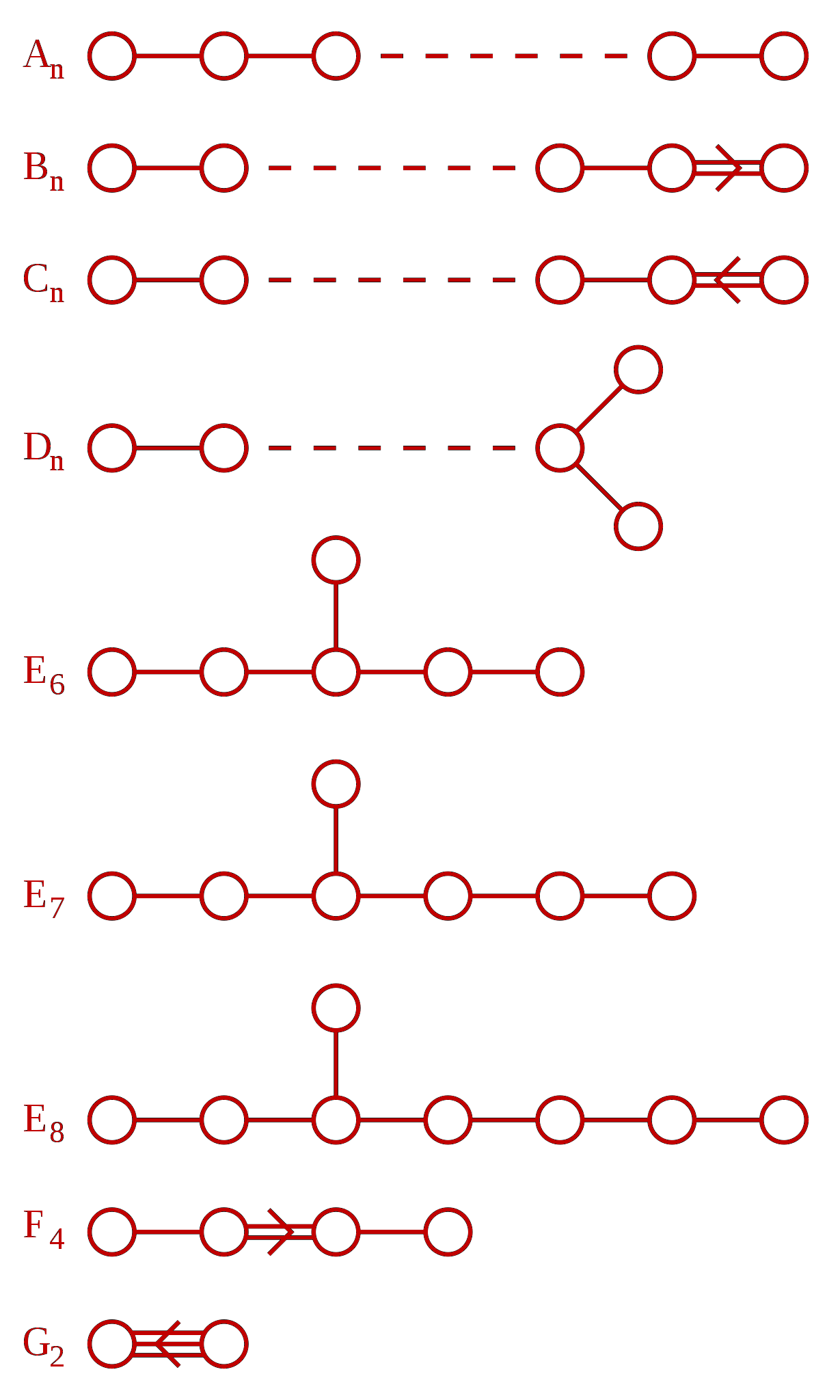} \\
   \tiny{\href{https://commons.wikimedia.org/wiki/File:Connected_Dynkin_Diagrams.svg}{image by R.\ A.\ Nonenmacher, CC BY-SA 4.0}}

  \end{center}
  
\newpage
\setcounter{page}{1}

I want to explain about a fascinating subject of huge importance in both
mathematics and physics: Coxeter and Dynkin diagrams.    Not to keep you in
suspense, I've already shown you the Dynkin diagrams on the front page of these notes.
Coxeter diagrams are a bit different, but similar. 

I'll start by explaining how Coxeter
diagrams can be used to classify ``finite reflection groups'': that is, finite groups
of linear transformations of \(\R^n\) generated by reflections.   Then I'll explain how 
Dynkin diagrams classify lattices in \(\R^n\) having finite reflection groups as symmetries.
Next, I'll explain how we use Dynkin diagrams to classify compact simple Lie groups and
their Lie algebras---and what those things are!    Finally, I'll say how
the ``simply-laced'' Dynkin diagrams---the ones without arrows on them---can be
used to classify integral lattices with a basis of vectors $v$ with $\|v\|^2 = 2$, quivers with tame representation theory, and finite subgroups of the group of rotations of 3-dimensional Euclidean space.

In short, Coxeter and Dynkin diagrams show up all over the place when you start trying to 
classify beautiful and symmetrical things.   In fact, we'll see that
the Platonic solids:
\begin{center}
\includegraphics[scale=0.1]{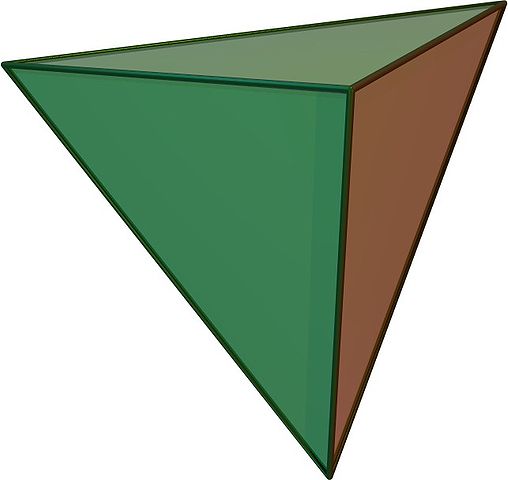} \vskip 0.7em
\includegraphics[scale=0.1]{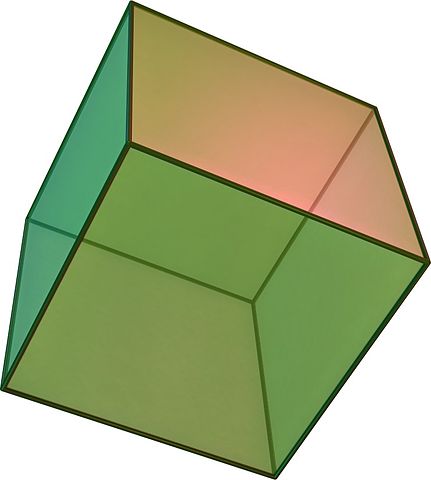}  \qquad \qquad \qquad \qquad \quad
\includegraphics[scale=0.1]{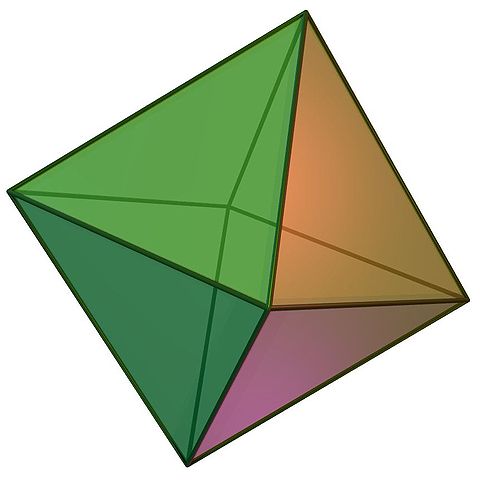}   \vskip 2.5em
\includegraphics[scale=0.1]{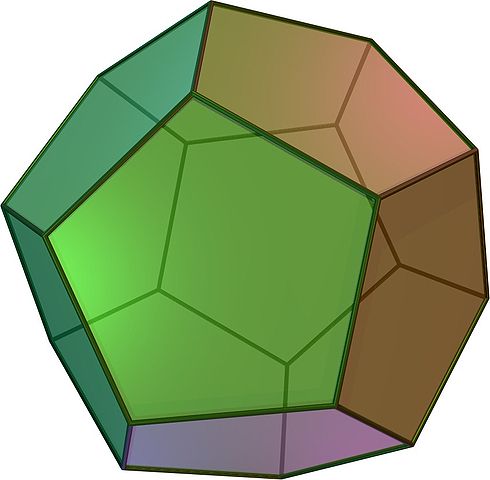}  \qquad \qquad 
\includegraphics[scale=0.1]{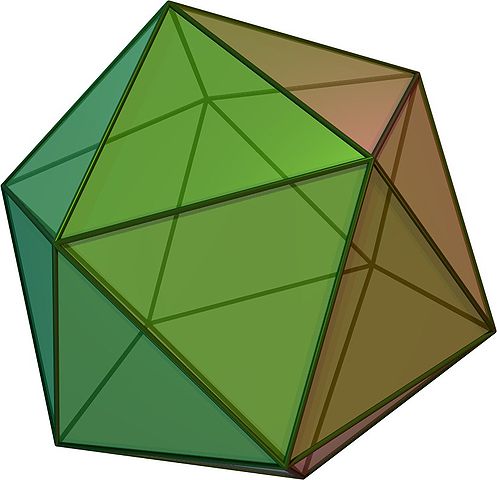}
\end{center}
are connected to Coxeter and Dynkin diagrams in \emph{two separate ways!}

\subsubsection*{Coxeter diagrams and finite reflection groups}

Okay, so what the heck is a Coxeter diagram?  We get these when we
try to classify ``finite reflection groups.'' Say
we are in \(n\)-dimensional Euclidean space. Then given any nonzero vector
\(v\), there is a reflection that takes \(v\) to \(-v\) and doesn't do
anything to the vectors orthogonal to \(v\). Let's call this a
\define{reflection through \(v\)}. A finite reflection group is a finite
group of transformations of Euclidean space such that every element is a
product of such reflections.  For example, the group of symmetries of a regular 
\(n\)-gon is a finite reflection group.  Showing this is a useful exercise if you 
don't see it right off the bat.

Note that if we do two reflections, we get a rotation. In particular,
suppose we have vectors \(v\) and \(w\) at an angle \(\theta\) from each
other, and let \(r\) and \(s\) be the reflections through \(v\) and
\(w\), respectively.  Then \(rs\) is a rotation by the angle \(2\theta\). Draw
a picture and check it! This means that if \(\theta = \pi /m\), then
\((rs)^n\) is a rotation by the angle \(2\pi\), which is the same as no
rotation at all, so \((rs)^m = 1\). On the other hand, if \(\theta\) is not a
rational number times \(\pi\), we never have \((rs)^n = 1\), so \(r\)
and \(s\) can not both be in some \emph{finite} reflection group.

So, if \(r\) and \(s\) are two distinct reflections in a finite reflection group 
we must have 
\[           r^2 = s^2 = 1\]
since doing a reflection twice gets you back where you started, but also
\[          (rs)^m = 1 \]
for some \(m = 2, 3, \dots \,\).

Using this, we can see that any finite
reflection group has a presentations described by a \define{Coxeter diagram}. 
The idea is that the group is generated by reflections through vectors that are 
at angles of \(\pi/m\) from each other for various choices of \(m = 2,3,
\dots\,\). To keep track of this, we draw
a dot for each one of these vectors.   Then, suppose two of the vectors are at
an angle \(\pi/m\) from each other.   If \(m = 2\), the reflections must commute,
and we don't bother drawing a line between the two dots.  Otherwise we draw a line between
them and label it with the number \(m\). 

Conversely, if someone hands us a Coxeter diagram
we get a group called a \define{Coxeter group} with 
\begin{itemize}
\item one generator \(r\) for each dot,
\item one relation \(r^2 = 1\) for each generator,
\item one relation \((rs)^2 = 1\) for each pair of generators with no line connecting their dots,
\item one relation \((rs)^m = 1\) for each pair of generators with a 
line labeled by the number \(m\) connecting their dots.
\end{itemize}

For example, the Coxeter diagram
 \[
  \begin{tikzpicture}
    \draw[very thick] (0,0) node{$\sbullet$} to (1,0) node{$\sbullet$} to node[label=above:{7}]{} (2,0) node {$\sbullet$};
    \draw[very thick] (0,0) to (0.5,-0.866) node{$\sbullet$} to (1,0);
        \node at (0, -0.6) {$4$};
        \node at (0.5, 0.4) {$3$};
        \node at (1, -0.6) {$3$};
  \end{tikzpicture}
\] 
gives us the Coxeter group with presentation
\[  \begin{array}{c} r^2 = s^2 = t^2 = u^2 = 1 \\
  (rs)^3 = (rt)^2 =  (ru)^4 = (st)^7  = (su)^3 = (tu)^2 = 1  \\
\end{array} \]
However, in this game a lot of edges wind up being labeled with the number \(3\).  People
usually leave out the label when this happens.    Then we draw the above diagram like this:
 \[
  \begin{tikzpicture}
    \draw[very thick] (0,0) node{$\sbullet$} to (1,0) node{$\sbullet$} to node[label=above:{7}]{} (2,0) node {$\sbullet$};
    \draw[very thick] (0,0) to (0.5,-0.866) node{$\sbullet$} to (1,0);
        \node at (0, -0.6) {$4$};
  \end{tikzpicture}
\] 

Now for some big theorems.  First, it turns out that if a Coxeter group is \emph{finite}, 
it's a finite reflection group!
Not every diagram yields a finite group.   But we can classify all 
possible Coxeter diagrams giving finite groups!  They have names, and they are
famous: they're like the rock stars of finite group theory.

First there is \(\mathrm{A}_n\), which has \(n\) dots in a row, like this: 
\[
  \begin{tikzpicture}
    \draw[very thick] (0,0) node{$\sbullet$} to (1,0) node{$\sbullet$} to (2,0) node{$\sbullet$} to (3,0) node{$\sbullet$};
  \end{tikzpicture}
\] 
For example, the group of symmetries of the equilateral triangle is the Coxeter
group \(\mathrm{A}_2\).  The group of symmetries of a regular tetrahedron
is \(\mathrm{A}_3\).
\begin{center}
\includegraphics[scale=0.2]{tetrahedron.jpg} 
\end{center}
More generally, \(\mathrm{A}_n\) is the group of symmetries of a regular
\(n\)-dimensional simplex ---which is just the group of all permutations of
the \(n+1\) vertices.   

Then there is \(\mathrm{BC}_n\), which has \(n\) dots:
\[
  \begin{tikzpicture}
    \draw[very thick] (0,0) node{$\sbullet$} to (1,0) node{$\sbullet$} to (2,0) node{$\sbullet$} to node[label=above:{4}]{} (3,0) node{$\sbullet$};
  \end{tikzpicture}
\] 
Only the last edge is labeled with a 4.   Later we'll see how  \(\mathrm{BC}_n\) spawns 
two Dynkin diagrams called  \(\mathrm{B}_n\) and  \(\mathrm{C}_n\); this accounts for
its strange name.   

The Coxeter group \(\mathrm{BC}_2\) is the
symmetry group of a square, while \(\mathrm{BC}_3\) is the symmetry group of the cube
or its dual, the regular octahedron.
\begin{center}
\includegraphics[scale=0.2]{hexahedron.jpg}  \qquad \qquad
\includegraphics[scale=0.2]{octahedron.jpg} 
\end{center}
In general, \(\mathrm{BC}_n\) is the group of symmetries of an \(n\)-dimensional
hypercube, or its dual, whose vertices lie at the centers of the faces of this hypercube.  
The dual of a hypercube could be called a ``regular 
hyperoctahedron'', but it's actually called a \define{cross-polytope} or \define{orthoplex}.
Here's how the duality works for \(n = 3\): 
\begin{center}
\includegraphics[scale=0.4]{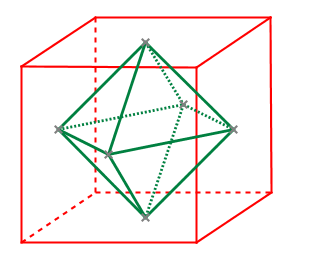}  
\end{center}

The next infinite series of finite reflection groups is \(\mathrm{D}_n\), which has \(n\) dots arranged in a row that branches at the end.  Here is \(\mathrm{D}_6\):
\[
  \begin{tikzpicture}
    \draw[very thick] (0,0) node{$\sbullet$} to (1,0) node{$\sbullet$} to (2,0) node{$\sbullet$} to (3,0) node{$\sbullet$};
    \draw[very thick] (3,0) to (3.87,0.87) node{$\sbullet$};
    \draw[very thick] (3,0) to (3.87,-0.87) node{$\sbullet$};
  \end{tikzpicture}
\] 
Since \(\mathrm{D}_n \cong \mathrm{A}_n\)
for \(n < 4\), the first really exciting case is \(\mathrm{D}_4\):
\[
  \begin{tikzpicture}
    \draw[very thick]  (2,0) node{$\sbullet$} to (3,0) node{$\sbullet$};
    \draw[very thick] (3,0) to (3.87,0.87) node{$\sbullet$};
    \draw[very thick] (3,0) to (3.87,-0.87) node{$\sbullet$};
  \end{tikzpicture}
\] 
and the symmetry of this diagram, called \define{triality}, gives birth to many remarkable
things. 

The Coxeter group \(\mathrm{D}_n\) is the symmetry group of the 
\(n\)-dimensional \define{demihypercube}, whose vertices are gotten by taking every other vertex of  a hypercube: that is, half the hypercube's vertices, with no two right next to
each other.  This is easiest to visualize when \(n = 3\).   
Here are the two demicubes in a cube:
\begin{center}
\includegraphics[scale=0.1]{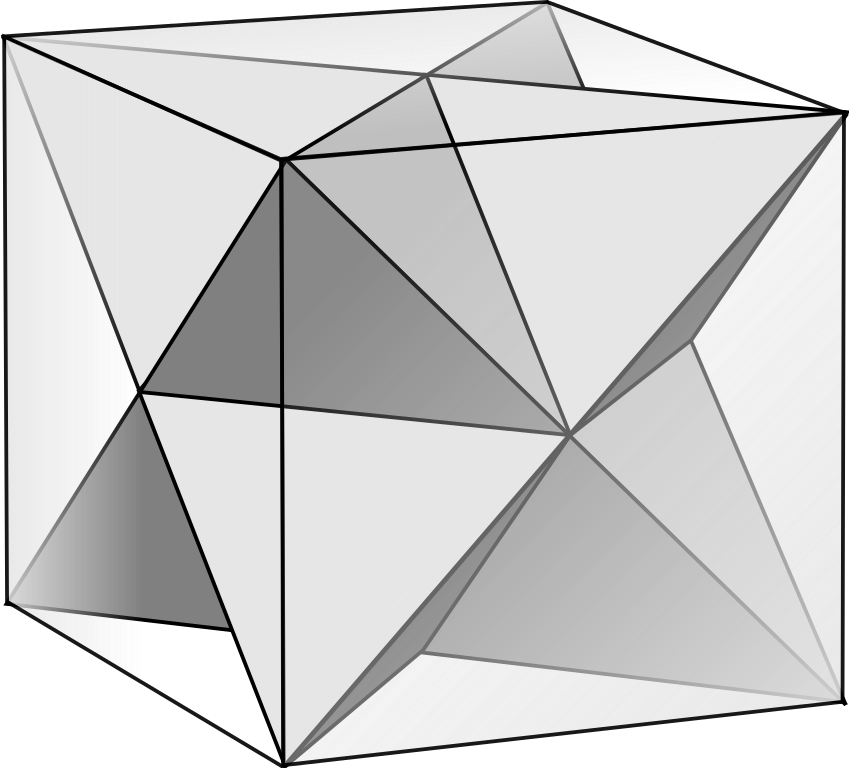}
\end{center}
They're just regular tetrahedra!   We might have expected this, 
since \(\mathrm{D}_3 \cong \mathrm{A}_3\).  A 4-dimensional demihypercube
has one of these tetrahedra ``on top'' and the other ``on the bottom''.   

Then there are \(\mathrm{E}_6\), \(\mathrm{E}_7\), and
\(\mathrm{E}_8\): 
\[
  \begin{gathered}
    \begin{tikzpicture}
      \draw[very thick] (0,0) node{$\sbullet$} to (1,0) node{$\sbullet$} to (2,0) node{$\sbullet$} to (3,0) node{$\sbullet$} to (4,0) node{$\sbullet$};
      \draw[very thick] (2,0) to (2,-1) node{$\sbullet$};
    \end{tikzpicture}
\\\begin{tikzpicture}
      \draw[very thick] (0,0) node{$\sbullet$} to (1,0) node{$\sbullet$} to (2,0) node{$\sbullet$} to (3,0) node{$\sbullet$} to (4,0) node{$\sbullet$} to (5,0) node{$\sbullet$};
      \draw[very thick] (3,0) to (3,-1) node{$\sbullet$};
    \end{tikzpicture}
  \\\begin{tikzpicture}
      \draw[very thick] (0,0) node{$\sbullet$} to (1,0) node{$\sbullet$} to (2,0) node{$\sbullet$} to (3,0) node{$\sbullet$} to (4,0) node{$\sbullet$} to (5,0) node{$\sbullet$} to (6,0) node{$\sbullet$};
      \draw[very thick] (4,0) to (4,-1) node{$\sbullet$};
    \end{tikzpicture}
  \end{gathered}
\] 
Interestingly, this series does \emph{not} go on.  Thus the \(\mathrm{A}_n\),
\(\mathrm{B}_n\), and \(\mathrm{D}_n\) Coxeter diagrams are called ``classical''
while the \(\mathrm{E}\) series is called ``exceptional''.   

A \define{polytope} is the higher-dimensional generalization of a polygon or polyhedron.
The \(\mathrm{E}_6,\mathrm{E}_7\) and \(\mathrm{E}_8\) Coxeter groups
are symmetries of some remarkably subtle polytopes:
\begin{itemize}
\item The \(\mathrm{E}_8\) Coxeter group is the 
symmetry group of an 8-dimensional polytope with 240 vertices 
called the \define{\(\mathrm{E}_8\) root polytope}.   To build this, take a 
sphere 8 dimensions and get as many equal-sized  
spheres as possible to touch it.   There will be 240.  The centers of these
spheres are the vertices of an \(\mathrm{E}_8\) root polytope.
\item The \(\mathrm{E}_7\) Coxeter group is the 
symmetry group of a 7-dimensional polytope with 126 vertices called the
\define{\(\mathrm{E}_7\) root polytope}.   To build this, pick 
any vertex of the \(\mathrm{E}_8\) root 
polytope.  Draw a line through it and the center of that polytope.  The vertices in the 7d
space orthogonal to this line are the vertices of a \(\mathrm{E}_7\) root polytope. 
\item The \(\mathrm{E}_6\) Coxeter group is the 
symmetry group of a 6-dimensional polytope with 72 vertices called the
\define{\(\mathrm{E}_6\) root polytope}.   Pick any 6 vertices of the \(\mathrm{E}_8\) root 
polytope forming a regular hexagon in some plane containing the center of that polytope. 
The vertices in the 6d space orthogonal to this plane are
the vertices of an \(\mathrm{E}_7\) root polytope. 
\end{itemize}
We will eventually give more manageable descriptions of these polytopes.

Then there is \(\mathrm{F}_4\): 
\[
  \begin{tikzpicture}
    \draw[very thick] (0,0) node{$\sbullet$} to (1,0) node{$\sbullet$} to node[label=above:{4}]{} (2,0) node{$\sbullet$} to (3,0) node{$\sbullet$};
  \end{tikzpicture}
\] 
The \(\mathrm{F}_4\)  Coxeter group is the symmetry group of a 4-dimensional 
polytope with 24 vertices and 24 octahedral faces, called the \define{24-cell}.   To build a
24-cell, first draw a 4-dimensional hypercube centered at the origin.   Put a point in the center of each face to get the vertices of a 4d orthoplex.   Then expand the orthoplex until its vertices
are as far from origin as the hypercube's vertices.   Then the vertices of the hypercube and orthoplex, taken together, are the vertices of a 24-cell!   

Regular polytopes are the most obvious higher-dimensional generalization of Platonic solids.
In all dimensions \(n > 4\) there are just three regular polytopes: the \(n\)-simplex, the \(n\)-dimensional hypercube and the \(n\)-dimensional orthoplex.  In dimension 4 there are three more!   One is the 24-cell.   We shall meet the other two soon.

Then there is \(\mathrm{G}_2\): 
\[
  \begin{tikzpicture}
    \draw[very thick] (0,0) node{$\sbullet$} to node[label=above:{6}]{} (1,0) node{$\sbullet$};
  \end{tikzpicture}
\] 
The \(\mathrm{G}_2\) Coxeter group is just the symmetry group of the regular hexagon.
There are infinitely many regular polygons; we'll see later why the equilateral triangle
(with symmetry group \(\mathrm{A}_2\), the square (with symmetry group \(\mathrm{B}_2\)
and the hexagon (with symmetry group \(\mathrm{G}_2\) are singled out as special.  If
you've ever tried to tile your floor with regular pentagons or heptagons maybe you can guess!

Then there are
\(\mathrm{H}_3\) and \(\mathrm{H}_4\): 
\[
\begin{tikzpicture}
    \draw[very thick] (0,0) node{$\sbullet$} to node[label=above:{5}]{} (1,0) node{$\sbullet$} to (2,0) node{$\sbullet$};
  \end{tikzpicture}
  \qquad
  \begin{tikzpicture}
  \draw[very thick] (0,0) node{$\sbullet$} to node[label=above:{5}]{} (1,0) node{$\sbullet$} to (2,0) node{$\sbullet$} to (3,0) node{$\sbullet$};
\end{tikzpicture}
\] 
\(\mathrm{H}_3\) is the group of symmetries of the regular dodecahedron or
icosahedron. 
\begin{center}
\includegraphics[scale=0.2]{dodecahedron.jpg}  \qquad \qquad
\includegraphics[scale=0.2]{icosahedron.jpg} 
\end{center}
\(\mathrm{H}_4\) is the group of symmetries of the remaining two 4d regular polytopes.  
They could be called the ``hyperdodecahedron'' and ``hypericosahedron'', 
but in fact they are called the 120-cell and 600-cell.   

To get your hands on these things, it is quickest to start with a regular 
dodecahedron.   It has 60 rotational symmetries since you can rotate your
favorite faces to any of the 12 faces, and in 5 ways.
These rotational symmetries form a group \(\Gamma\) that's a subgroup of \(\SO(3)\), the 
group of all rotations in 3d space.   There's a 2-1 map from \(\SU(2)\) to \(\SO(3)\), so two 
elements of \(\SU(2)\) map to each element of \(\Gamma\).  These element of 
\(\SU(2)\) form a 120-element subgroup of \(\SU(2)\).  But geometrically \(\SU(2)\) is a 
sphere in 4d space: for example, it's isomorphic to the group of unit quaternions.   So, 
\(\Gamma\) consists of 120 points on a sphere in 4d space... and these are the
vertices of the \define{600-cell}.

The 600-cell has 120 vertices and 600 tetrahedral faces.  If we put a 
point at the center of each face, we get the vertices of the dual polytope, which has
600 vertices and 120 tetrahedral faces.   This is called the \define{120-cell}.


Finally, there is another very simple infinite series, \(\mathrm{I}_m\): 
\[
  \begin{tikzpicture}
    \draw[very thick] (0,0) node{$\sbullet$} to node[label=above:{$m$}]{} (1,0) node{$\sbullet$};
  \end{tikzpicture}
\] 
The Coxeter group coming from this diagram is the symmetry group of the
regular \(m\)-gon.  Notice that
\[    \mathrm{I}_3 \cong \mathrm{A}_2, \quad  
\mathrm{I}_4 \cong \mathrm{B}_2, \quad
\mathrm{I}_6 \cong \mathrm{G}_2.  \]
So, the only \emph{new} Coxeter diagrams here are \(\mathrm{I}_m\) for \(m = 5\) 
and \(m \ge 7\).

We can get other finite reflection groups from disjoint unions of the Coxeter
diagrams we've already seen: they're just products of the Coxeter groups we've
already seen.  But our list of \emph{connected} Coxeter diagrams giving finite reflection
groups is \emph{done!}  

Let's list them without repetitions:
\begin{itemize}
\item
  \(\mathrm{A}_n\), \(n > 0\)
  \[
  \begin{tikzpicture}
    \draw[very thick] (0,0) node{$\sbullet$} to (1,0) node{$\sbullet$} to (2,0) node{$\sbullet$} to (3,0) node{$\sbullet$};
  \end{tikzpicture}
\] 
\item
  \(\mathrm{BC}_n\), \(n > 1\)
  \[
  \begin{tikzpicture}
    \draw[very thick] (0,0) node{$\sbullet$} to (1,0) node{$\sbullet$} to (2,0) node{$\sbullet$} to node[label=above:{4}]{} (3,0) node{$\sbullet$};
  \end{tikzpicture}
\] 
\item
  \(\mathrm{D}_n\), \(n > 3\)
  \[
  \begin{tikzpicture}
    \draw[very thick] (0,0) node{$\sbullet$} to (1,0) node{$\sbullet$} to (2,0) node{$\sbullet$} to (3,0) node{$\sbullet$};
    \draw[very thick] (3,0) to (3.87,0.87) node{$\sbullet$};
    \draw[very thick] (3,0) to (3.87,-0.87) node{$\sbullet$};
  \end{tikzpicture}
\] 
\item
  \(\mathrm{E}_6\), \(\mathrm{E}_7\), \(\mathrm{E}_8\)
  \[
  \begin{gathered}
    \begin{tikzpicture}
      \draw[very thick] (0,0) node{$\sbullet$} to (1,0) node{$\sbullet$} to (2,0) node{$\sbullet$} to (3,0) node{$\sbullet$} to (4,0) node{$\sbullet$};
      \draw[very thick] (2,0) to (2,-1) node{$\sbullet$};
    \end{tikzpicture}
\\\begin{tikzpicture}
      \draw[very thick] (0,0) node{$\sbullet$} to (1,0) node{$\sbullet$} to (2,0) node{$\sbullet$} to (3,0) node{$\sbullet$} to (4,0) node{$\sbullet$} to (5,0) node{$\sbullet$};
      \draw[very thick] (3,0) to (3,-1) node{$\sbullet$};
    \end{tikzpicture}
  \\\begin{tikzpicture}
      \draw[very thick] (0,0) node{$\sbullet$} to (1,0) node{$\sbullet$} to (2,0) node{$\sbullet$} to (3,0) node{$\sbullet$} to (4,0) node{$\sbullet$} to (5,0) node{$\sbullet$} to (6,0) node{$\sbullet$};
      \draw[very thick] (4,0) to (4,-1) node{$\sbullet$};
    \end{tikzpicture}
  \end{gathered}
\] 
\item
  \(\mathrm{F}_4\)
\[
  \begin{tikzpicture}
    \draw[very thick] (0,0) node{$\sbullet$} to (1,0) node{$\sbullet$} to node[label=above:{4}]{} (2,0) node{$\sbullet$} to (3,0) node{$\sbullet$};
  \end{tikzpicture}
\] 
\item
  \(\mathrm{G}_2\)
  \[
  \begin{tikzpicture}
    \draw[very thick] (0,0) node{$\sbullet$} to node[label=above:{6}]{} (1,0) node{$\sbullet$};
  \end{tikzpicture}
\] 
\item
  \(\mathrm{H}_4\), \(\mathrm{H}_5\)
\[
\begin{tikzpicture}
    \draw[very thick] (0,0) node{$\sbullet$} to node[label=above:{5}]{} (1,0) node{$\sbullet$} to (2,0) node{$\sbullet$};
  \end{tikzpicture}
  \qquad
  \begin{tikzpicture}
  \draw[very thick] (0,0) node{$\sbullet$} to node[label=above:{5}]{} (1,0) node{$\sbullet$} to (2,0) node{$\sbullet$} to (3,0) node{$\sbullet$};
\end{tikzpicture}
\] 
\item
 \(\mathrm{I}_m\)
 \[
  \begin{tikzpicture}
    \draw[very thick] (0,0) node{$\sbullet$} to node[label=above:{$m$}]{} (1,0) node{$\sbullet$};
  \end{tikzpicture}
\] 
\end{itemize}
The proof that these are the only possibilities is a rather elaborate inductive argument.  For details, see for example \textsl{Finite Reflection Groups and Coxeter Groups} by Humphreys.

Note that the symmetry groups of the Platonic
solids and their higher-dimensional relatives fit in nicely into
this classification.   Later we shall see a more mysterious relation between
Platonic solids and these diagrams.  But first we need to introduce Dynkin 
diagrams!  These show up when we think about lattices.

\subsubsection*{Dynkin diagrams and lattices}

We get a \define{lattice} by taking \(n\) linearly independent vectors in \(n\)-dimensional Euclidean space and forming all linear combinations with integer coefficients.   
\begin{center}
\includegraphics[scale=0.3]{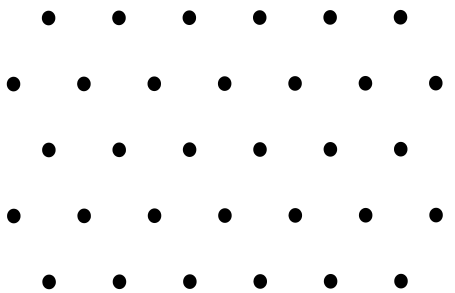}
\end{center}
Sometimes lattices have interesting symmetry groups.  Every lattice in \(n\) dimensions has 
\(\mathbb{Z}^n\) acting as translation symmetries, so let's focus on symmetries that
fix the origin: that is, combinations of rotations and reflections that map the lattice
to itself.    To do this, we can exploit the classification of finite reflection groups.

So, suppose we have a connected Coxeter diagram that gives a finite reflection
group.   If our diagram has \(n\) dots, this group acts on \(\R^n\).   \emph{When is there a
lattice in \(\R^n\) having this group as symmetries?}     

If one exists, we say our group satisfies the \define{crystallographic
condition}.    It turns out that the only ones that do are
\[\mbox{$\mathrm{A}_n$, \(\mathrm{B}_n$, \(\mathrm{D}_n$, \(\mathrm{E}_6$, \(\mathrm{E}_7$, \(\mathrm{E}_8$, \(\mathrm{F}_4$, and \(\mathrm{G}_2$}!\]
In other words, Coxeter diagrams with any edges labeled by numbers \(m = 5\) or \(m \ge 7\)
are ruled out.  So, we can't find lattices whose symmetries include the symmetry groups of the
regular pentagon, the regular heptagon, or regular polygons with more sides.   This has big implications for crystals, but also for pure math.

Say we have a finite reflection group \(\Gamma\) acting on \(\R^n\) and it obeys the crystallographic condition.  How can get a lattice in \(\R^n\) with \(\Gamma\) as symmetries?   It turns out we can always do it by taking integer linear combinations of some basis vectors, one for each dot in the Coxeter diagram.  The reflections generating \(\Gamma\) will be reflections through these vectors.

To get this to work, the angle between two of these vectors needs to be \(\pi/m\) when the edge between the dots is labeled by \(m\).    But remember: in this game an unlabeled edge counts 
as an edge labeled by \(3\), so then the vectors need to be at a \(\pi/3\) angle from each 
other.  No edge at all counts as an edge labeled by \(2\), so then the vectors need to be
at an angle of \(\pi/2\) from each other---that is, at right angles.

For example, take \(\mathrm{A}_2\):  
  \[
  \begin{tikzpicture}
    \draw[very thick] (0,0) node{$\sbullet$}  to node[label=above:{3}]{} (1,0) node{$\sbullet$};
  \end{tikzpicture}
\]
where I'm writing the \(3\) to clarify some patterns.
To get a lattice with this Coxeter group as symmetries, take two vectors of the
same length at an angle of \(\pi/3\) from each other.   Their integer linear
combinations form the lattice we want.  It's the lattice of vertices
of the \define{triangular tiling}:
\begin{center}
\includegraphics[scale=0.3]{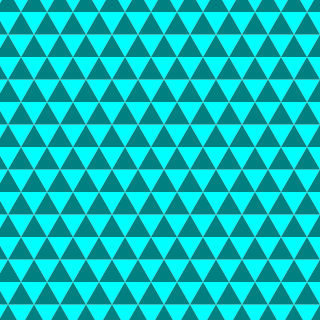}
\end{center}

Or consider \(\mathrm{B}_2\):  
\[
  \begin{tikzpicture}
    \draw[very thick] (0,0) node{$\sbullet$} to node[label=above:{4}]{} (1,0) node{$\sbullet$};
  \end{tikzpicture}
\]
To get a lattice with this Coxeter group as symmetries, we can use
two vectors at an angle of \(\pi/4\) from each other.   But they can't be of
equal length!  One must be \(\sqrt{2}\) times as long as the other.   Then
their integer linear combinations form the lattice we want.  
It's the lattice of vertices of the \define{square tiling}:
\begin{center}
\includegraphics[scale=0.15]{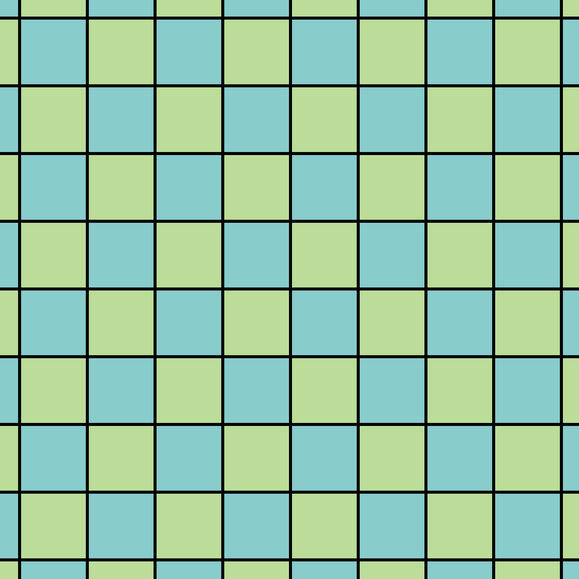}
\end{center}

Similarly, for \(\mathrm{G}_2\):  
\[
  \begin{tikzpicture}
    \draw[very thick] (0,0) node{$\sbullet$} to node[label=above:{6}]{} (1,0) node{$\sbullet$};
  \end{tikzpicture}
\]
To get a lattice with this Coxeter group as symmetries, we can use
two vectors at an angle of \(\pi/6\) from each other.  One must be \(\sqrt{3}\) 
times as long as the other.   Then their integer linear combinations form a lattice with the symmetries we want.   It's the lattice of vertices of the \define{hexagonal tiling}:
\begin{center}
\includegraphics[scale=0.3]{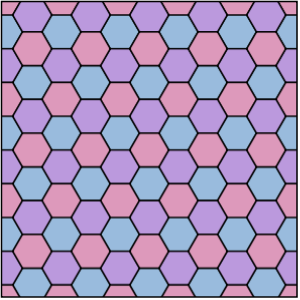}
\end{center}
No, it's not!  I was just checking to see if you're paying attention!  The vertices of
the hexagonal tiling don't form a lattice.  To get a lattice you can throw in an
extra point at the center of each hexagon, but then you get the triangular tiling
again.   In fact the \(\mathrm{G}_2\) lattice looks just like the \(\mathrm{A}_2\) lattice, 
at least up to a rotation and rescaling.  The difference is that now we're considering it 
with a bigger symmetry group: \(\mathrm{A}_2\) gives the symmetry group of an equilateral 
triangle, while \(\mathrm{G}_2\) gives the symmetry group of a regular hexagon, which is twice
as big.

Now suppose we have a connected Coxeter diagram that gives a finite reflection group \(\Gamma\).   And suppose
it obeys the crystallographic conditions, so all  its edges are labeled by numbers \(m = 3,4\) or \(6\). 
Here is how we build a lattice called the \define{root lattice} having \(\Gamma\) as its symmetries.

Suppose our Coxeter diagram has \(n\) dots.  Then we can pick a basis of vectors \(v_1, \dots, v_n \in \R^n\), one for each dot in the diagram, with these properties:
\begin{enumerate}[label=(\alph*)]
\item If there is no edge connecting two dots, the vectors for these dots are at an angle \(\pi - \pi/2\) from each other.
\item If there is an edge with \(m = 3\) connecting two dots, their vectors are at an angle \(\pi - \pi/3\) from each other, and have the same length.
\item If there is an edge with \(m = 4\) connecting two dots, their vectors are at an angle \(\pi - \pi/4\) from each other, and  one vector is\(\sqrt{2}\) times as long as the other. 
\item If there is an edge with \(m = 6\) connecting two dots, their vectors are at an angle \(\pi - \pi/6\) from each other, and one vector is \(\sqrt{3}\) times as long as the other.
\end{enumerate}
You're probably wondering why we are choosing vectors that are angles \(\pi - \theta\) from each other, when you were expecting the angle to be \(\theta\).   Of course if the angle between vectors $v$ and $w$ is $\pi - \theta$, the angle between $v$ and $-w$ is $\theta$.  So why are we choosing to work with $\pi  - \theta$ above?  This is the convention everyone uses.  Understanding why this convention is better requires a deeper study of ``root systems'', which I am trying to sidestep here.

The recipe above is missing an ingredient: when \(m = 4\) or \(6\) we need to decide which vector is longer.  To make this decision, let us draw an arrow on the edge from the dot with the longer vector to the dot with the shorter vector.   These arrows are enough to specify a lattice on which our Coxeter group acts as symmetries!---at least up to rotations, reflections and rescalings.

For example, in the diagram \(\mathrm{F}_4\) we have two choices:
\[
  \begin{tikzpicture}
    \draw[very thick] (0,0) node{$\sbullet$} to (1,0) node{$\sbullet$} to node[label=above:{4}]{\textgreater} (2,0) node{$\sbullet$} to (3,0) node{$\sbullet$};
  \end{tikzpicture}
\] 
and
\[
  \begin{tikzpicture}
    \draw[very thick] (0,0) node{$\sbullet$} to (1,0) node{$\sbullet$} to node[label=above:{4}]{\textless} (2,0) node{$\sbullet$} to (3,0) node{$\sbullet$};
  \end{tikzpicture}
\] 
However, in this case the resulting diagrams with arrows are isomorphic: you can turn
one around and get the other.

In fact, you can check that for every type of Coxeter diagram except \(\mathrm{BC}_n\) with \(n \ge 3\), you get isomorphic diagrams with arrows no matter how you point the arrows!   So, except in those cases, there is only \emph{one} lattice having that group of symmetries---up to rotations, reflections and rescalings.   But for \(\mathrm{BC}_n\) with \(n \ge 3\) there are two.  Recall that the diagram \(\mathrm{BC}_n\) looks like this, with \(n\) dots:
\[
  \begin{tikzpicture}
    \draw[very thick] (0,0) node{$\sbullet$} to (1,0) node{$\sbullet$} to (2,0) node{$\sbullet$} to node[label=above:{4}]{} (3,0) node{$\sbullet$};
  \end{tikzpicture}
\] with \(n\) dots.    
The two ways to draw arrows on the last edge are called \(\mathrm{B}_n\) and \(\mathrm{C}_n\):
\begin{itemize}
\item
  \(\mathrm{B}_n\): \[
      \begin{tikzpicture}
        \draw[very thick] (0,0) node{$\sbullet$} to (1,0) node{$\sbullet$} to (2,0) node{$\sbullet$} to node[label=above:{4}]{\textgreater} (3,0) node{$\sbullet$};
      \end{tikzpicture}
    \]
\item
  \(\mathrm{C}_n\): \[
      \begin{tikzpicture}
        \draw[very thick] (0,0) node{$\sbullet$} to (1,0) node{$\sbullet$} to (2,0) node{$\sbullet$} to node[label=above:{4}]{\textless} (3,0) node{$\sbullet$};
      \end{tikzpicture}
    \]
\end{itemize}
The points near the origin in the \(\mathrm{B}_3\) lattice lie on a cube,
while those in the \(\mathrm{C}_3\) lie on an octahedron:

\begin{center}
\includegraphics[scale = 0.8]{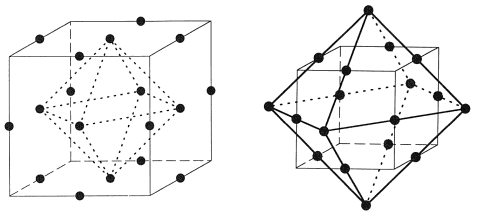}
\end{center}

\noindent
That's how it works in higher dimensions too, with a hypercube and orthoplex.

People also use another style of diagram to describe what's going on here: ``Dynkin diagrams''.
For these, we take our Coxeter diagram with arrows on it, and:
\begin{itemize}
\item replace any edge labeled with a 4 by two parallel edges;
\item replace any edge labeled with a 6 by three parallel edges.
\end{itemize}

Yes, I know this sounds confusing at first!   It's just a different way to draw the same
information.   I guess people like it because you can draw pictures without any numbers. Let me do a couple of examples.  Our friend \(\mathrm{G}_2\), with an arrow on it:
\[
  \begin{tikzpicture}
    \draw[very thick] (0,0) node{$\sbullet$} to node[label=above:{6}]{\textless} (1,0) node{$\sbullet$};
  \end{tikzpicture}
\]
gets drawn as this Dynkin diagram:
 \[
  \begin{tikzpicture}
    \draw[thick] (0,0) to (1,0);
    \draw[thick] (0,0.07) to (1,0.07);
    \draw[thick] (0,-0.07) to (1,-0.07);
    \node at (1, 0) {$\sbullet$};
        \node at (0, 0) {$\sbullet$};
     \node at (0.5, 0) {\scalebox{1.7}{$<$}};
  \end{tikzpicture}
\] 
while our friend \(B_3\):
\[
      \begin{tikzpicture}
        \draw[very thick] (1,0) node{$\sbullet$} to (2,0) node{$\sbullet$} to node[label=above:{4}]{\textgreater} (3,0) node{$\sbullet$};
      \end{tikzpicture}
    \]
gets drawn as this Dynkin diagram:
 \[
  \begin{tikzpicture}
  \draw[very thick] (-1,0) to (0,0);
    \draw[thick] (0,0.03) to (1,0.03);
    \draw[thick] (0,-0.03) to (1,-0.03);
    \node at (-1, 0) {$\sbullet$};
    \node at (1, 0) {$\sbullet$};
        \node at (0, 0) {$\sbullet$};
     \node at (0.5, 0) {\scalebox{1.7}{$>$}};
  \end{tikzpicture}
\] 

The upshot is that any Dynkin diagram with \(n\) dots describes a basis of
vectors \(v_1, \dots, v_n \in \R^n\) such that:
\begin{enumerate}
\item
integer linear combinations of these vectors form a lattice \(L \subset \R^n\)
\item 
reflections through these vectors generate a finite reflection group \(\Gamma\)
\item
the action of \(\Gamma\) on \(\R^n\) preserves the lattice \(L\).
\end{enumerate}
People usually normalize these vectors so that the shortest ones have \(v_i \cdot v_i = 2\), 
for reasons that should become clear later.   This determines them up to
rotations and reflections.   We then call these  vectors \(v_i\) \define{roots}, and we call the
lattice \(L\) a \define{root lattice}.  

If the action of \(\Gamma\) on \(\R^n\) is \define{indecomposable}, meaning we can't chop 
\(\R^n\) into a direct sum of two nontrivial subspaces both preserved by \(\Gamma\), then 
any collection of roots obeying 1--3 comes from a \emph{connected} Dynkin diagram!    And here 
are all the connected Dynkin diagrams, drawn in a more artistic style by R.\ A.\ Nonemacher.   The 
dots are drawn as big circles:

\begin{center}
 \includegraphics[scale = 0.14]{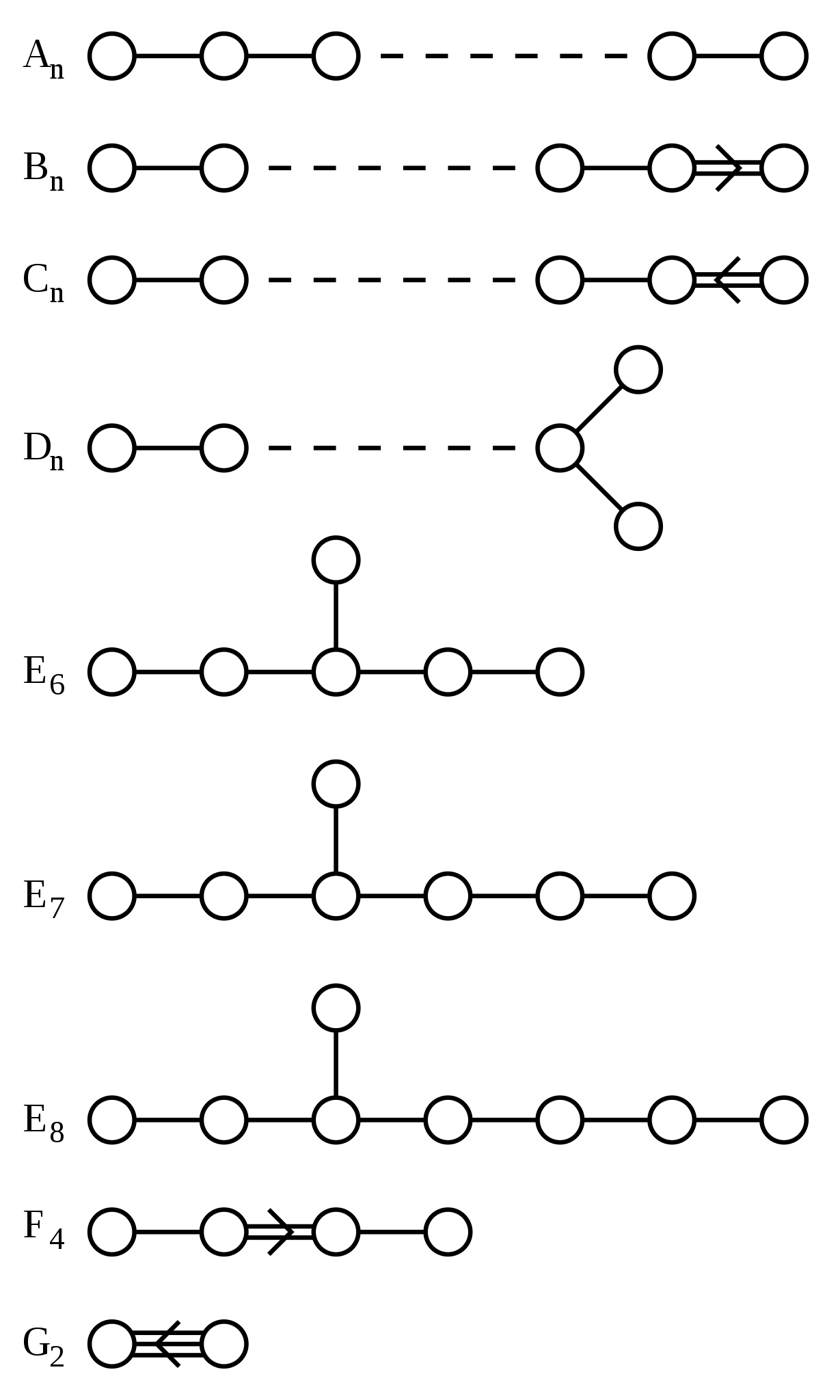}
\end{center}

\subsection*{Compact simple Lie algebras}

Now let us turn to the theory of Lie groups. Lie groups are the most
important ``continuous'' (as opposed to discrete) symmetry groups.
Examples include the real line with addition as the group operation,
the circle with addition modulo \(2\pi\), and the so-called ``classical
groups'', which include:
\begin{itemize}
\item
The \define{general linear group} \(\GL(N,\C)\), consisting of all invertible
linear transformations of \(\C^N$, or in other words, all \(N \times N\) complex
matrices with nonzero determinant.
\item 
The \define{special linear group} \(\SL(N,\C)\), consisting of all
linear transformations of \(\C^N\) with determinant \(1\).   
\item 
The \define{unitary group} \(\U(N)\), consisting of all unitary linear
transformations of \(\C^N\).
\item 
The \define{special unitary group} \(\SU(N)\), consisting of all unitary linear
transformations of \(\C^N\) with determinant \(1\).
\end{itemize}
All these Lie groups are incredibly
important in both physics and mathematics. Thus it is wonderful, and
charmingly ironic, that the same Dynkin diagrams that classify the
oh-so-discrete lattices with finite reflection groups as symmetries 
also classify some of the most
beautiful of Lie groups: the ``simple'' Lie groups.  

There is a vast amount known about semisimple Lie groups, and everyone
really serious about mathematics winds up needing to learn some of this
stuff. I took courses on Lie groups and their Lie algebras in grad
school, but it was only later that I really came to appreciate the
beauty of the simple Lie groups.  One reason I found them frustrating was that
the work involved in their classification was so algebraic, and I
preferred the more geometrical aspects of Lie groups.   This algebraic
approach de-emphasized the Lie groups themeselves, and emphasized
a tool for working with Lie group: namely, Lie \emph{algebras}.

So what's the basic idea? Let me summarize two semesters of grad school,
and tell you the basic stuff about Lie groups and the classification of
simple Lie groups. Forgive me if it's a bit rushed and sketchy: hopefully 
the main ideas will shine through the murk better this way.

A Lie group is a group that's also a manifold, for which the group
operations (multiplication and taking inverses) are smooth functions.
This lets you form the tangent space to any point in the group, and the
tangent space at the identity plays a special role. It's called the \define{Lie
algebra of the group}. If we have any element \(x\) in the Lie algebra,
we can exponentiate it to get an element \(\exp(x)\) in the group, and
we can keep track of the noncommutativity of the group by forming the
\define{Lie bracket} of elements \(x\) and \(y\) in the Lie algebra:
\[[x,y] = \left. \frac{d}{dt}\frac{d}{ds} \exp(sx) \exp(ty) \exp(-sx) \exp(-ty)\right|_{s,t = 0}\]
where \(s\) and \(t\) are real numbers.  Note that \([x,y] = 0\) if the group is 
commutative.
This bracket operation satisfies some axioms, and algebraists call
anything a Lie algebra that satisfies those axioms.  Just to satisfy your
curiosity, a \define{Lie algebra} is a vector space $\g$ with a bracket operation
such that
\begin{itemize}
\item the bracket is linear in the second argument:
\[    [x,ay + bz] = a[x,y] + b[x,z]   \qquad \textrm{ for all } x,y,z \in \g, \, a,b \in \R \]
\item the bracket is antisymmetric: 
\[   [x,y] = -[y,x] \qquad \textrm{ for all } x,y \in \g \]
and thus linear in the first argument as well,
\item bracketing with any element \(x\) satisfies a version of the product rule 
\[    [x,[y,z]] = [[x,y], z] + [y, [x,z]]  \qquad \textrm{ for all } x,y,z \in \g \]
called the \define{Jacobi identity}
\end{itemize}
All our Lie algebras will be finite-dimensional real vector spaces.  For example, you
could take \(n \times n\) real matrices and let \([x,y] = xy - yx\).  

I will not actually do anything with these axioms, since this is not really a course
on Lie algebras.  But these axioms are part of why Lie algebras are so wonderful.  
They're all about linear algebra---they're
just a vector space with a bracket operation obeying some axioms---and
yet we can \emph{almost} recover any Lie group from its Lie algebra.

Why ``almost''?   Well, first of all, we can't do it unless we require that our 
Lie group is connected.  The reason is that we could add extra connected components 
without changing the tangent space at the identity or its bracket operation.  
For example, every finite group is a Lie group with Lie algebra \(\{0\}\).   And 
second of all, we can't do it unless we require a bit
more.  The reason is that any covering space of a connected Lie group is 
another connected Lie group with the same Lie algebra.   Luckily, every connected 
Lie group has  a \define{universal cover}  which has no connected covering spaces except itself.  
This guy is both connected and simply connected.    

With these caveats we're okay: we can recover a connected and simply 
connected Lie group from its Lie algebra.  Even better, \emph{every} finite-dimensional
Lie algebra comes from a connected and simply connected Lie group!

This reduces the problem of classifying connected and 
simply connected Lie groups to a problem in linear algebra.  Unfortunately it's
an incredibly hard problem in linear algebra.   It appears to be hopeless unless
we stick to low dimensions.   There are just too many ways to build bigger Lie
groups, or Lie algebras, from smaller ones.   But the problem becomes much 
easier if we stick to \emph{compact} Lie groups---so that's what we will do.

The Lie algebra of a compact Lie group is called a \define{compact Lie algebra}.
It's not that the Lie algebra is literally compact: it's a vector space.  But 
compact Lie algebras are very nice.   It turns out that we can take direct sums of 
Lie algebras by defining operations componentwise, and a compact Lie algebra
is always the direct sum of an ``abelian'' Lie algebra and a ``semisimple'' one.
As we'll see, these lie at opposite extremes in a certain sense.

To understand this, it helps to think about the ``Killing form'' of a Lie algebra $\g$. 
For any $x \in \g$ there's a linear operator on $\g$ given by bracketing with $x$.
It's usually called $\ad_x \maps \g \to \g$, so
\[          \ad_x(y) = [x,y]  .\]
The \define{Killing form} is a bilinear form on \(\g\) given by
\[     B(x, y) = \tr(\ad(x) \ad(y))   .\]
By the cyclic property of the trace,
\[     B(x, y) = B(y , x).\]
This is great: it's a bit like getting an inner product for free!  But stay tuned: 
it's a bit subtler than that. 

Here are the two kinds of Lie algebras I mentioned:
\begin{itemize}
\item A Lie algebra is \define{abelian} if  \([x,y] = 0\) for all \(x\)
and \(y\).
\item A Lie algebra is \define{semisimple} if the Killing form is nondegenerate: if
\( \langle x,y \rangle = 0\) for all \(y\) then \(x = 0\).
\end{itemize}
Note that the Killing form of an abelian Lie algebra is zero, which is at the opposite
extreme from being nondegenerate.   The Killing form of a semisimple Lie algebra
is not really an inner product---it doesn't need to be positive definite.  But it comes very close 
for a \define{compact semisimple Lie algebra}---that is, a semisimple Lie algebra that's the Lie algebra of a compact Lie group.   Namely: for such a Lie algebra, the Killing form is
\emph{negative} definite, so 
\[    \langle x, y \rangle = -B(x,y) \]
is an inner product.  

There's one abelian Lie algebra of each dimension, and every abelian Lie algebra is the Lie
algebra of a \define{torus}, which is a product of finitely many circles.
So all the hard work lies in understanding the semisimple Lie algebras. 
A connected Lie group whose Lie algebra is semisimple is called---surprise!---a \define{semisimple Lie group}.   

Now let's start with a compact semisimple Lie group \(G\) and see how there's a Dynkin
diagram hiding inside it.    To do this, we'll construct a lattice in Euclidean space with a finite 
reflection group acting as symmetries.    Under favorable conditions this will be a root lattice, in 
the sense explained in the previous section.  

For starters, inside \(G\) there is always some subgroup \(T\) that's a torus not contained 
in any larger torus.  Such a subgroup is called a \define{maximal torus} of \(G\).   It's basically unique: 
there are a bunch of maximal tori, but any two  are conjugate to each  other in \(G\).  So, people 
often sloppily talk about  ``the''  maximal torus.

Let \(\mathrm{Lie}(T)\) stand for the Lie algebra of a maximal torus \(T\).
The inner product on \(\mathrm{Lie}(G)\) restricts to an inner product on \(\mathrm{Lie}(T)\).
So, \(\mathrm{Lie}(T)\) is isomorphic to good old \(\mathbb{R}^n\) with its
usual inner product.  Futhermore, sitting inside \(\mathrm{Lie}(T)\) there is a lattice \(L\), consisting of all elements  \(x\) with \(\exp(2 \pi x) = 1\).  This is called the \define{integer lattice} of \(G\).

Now, you might hope that this is the right way to build the desired root lattice.
Unfortunately it's not.  \(L\) is a very interesting lattice, but not quite
the one we want.  Luckily, it gives rise to a lattice in the dual 
vector space \(\mathrm{Lie}(T)^\ast\), namely
\[    \left\{ \ell \in \mathrm{Lie}(T)^\ast \, \vert \, \ell(v) \in \Z \textrm{ for all } v \in L \right\} .\]
This is called the \define{dual} of \(L\), and I'll denote it by \(L^\ast\).    

The dual lattice \(L^\ast\) is better for our purposes, but we're not quite out of the woods yet.   
As I already mentioned, any covering space of a connected Lie group is another 
connected Lie group with the same Lie algebra.    If we start with some compact semisimple 
Lie algebra \(\g\), we need to choose the right Lie group \(G\) with \(\mathrm{Lie}(G) = \g\) 
if we want to be sure that \(L^\ast\) is a root lattice.    One obvious guess is to let \(G\) be 
simply connected: this covers all the other connected Lie groups with Lie algebra \(\g\).  
However, if we do this, \(L^\ast\) is a lattice called the \define{weight lattice} of \(\g\).  It's very 
important, but it's not what we're after.   Instead, we should choose \(G\) so that its center is trivial: this is
\emph{covered by} all other connected Lie groups with Lie algebra \(G\).    Then \(L^\ast\)
is the root lattice we're after!

How do finite reflection groups get into the game?
For some elements \(g\) in \(G\), if we conjugate \(T\) by
\(g\), that is, form the set of all elements \(gtg^{-1}\) where \(t\) is
in \(T\), we get \(T\) back. In other words, these elements of \(G\) act
as symmetries of the torus \(T\).   But some elements of \(g \in G\)
act trivially on \(T\): they have \(gtg^{-1} = t\) for all \(t \in T\).   So, the 
quotient group
\[    W(G) = \frac{ \left\{g \in G \, \vert \, gtg^{-1} \in T \textrm{ for all } t \in T \right\} }
{ \left\{g \in G \, \vert \, gtg^{-1} = t  \textrm{ for all } t \in T \right\} } \]
acts on our maximal torus \(T\).  This group \(W(G)\) is called the \define{Weyl
group} of \(G\).   (In reality it also depends on our choice of maximal torus,
but changing our maximal torus to another one gives an isomorphic Weyl
group, so our notation ignores the dependence on \(T\).)

If a group acts as symmetries
of something, it also acts as symmetries of everything naturally
cooked up from that thing.  For this reason, the
Weyl group of \(G\) also acts as symmetries of \(\mathrm{Lie}(T)\), and the integer lattice $L$,
and the root lattice \(L^\ast\) sitting inside \(\mathrm{Lie}(T)^\ast\).  So we get a
lattice $L^\ast$, which we can think of as living in \(n\)-dimensional Euclidean space, 
together with a group of symmetries \(W(G)\).  

This has been a bit technical, and I apologize for that.  But now we are rewarded by 
three wonderful  theorems.  First, \(W(G)\) 
is actually a finite reflection group acting on the inner product space \(\textrm{Lie}(T)^\ast\).   Second, 
\(W(G)\) is generated by reflections through certain vectors in the lattice \(L^\ast\).
From what we learned in the previous section, 
this means the root lattice \(L^\ast\) and this finite reflection group \(W(G)\) are
determined, up to isomorphism, by some Dynkin diagram.   And third, it turns 
out that the Lie algebra of \(G\) is determined, up to isomorphism, by the root lattice \(L^\ast\) 
and this finite reflection group acting on it!

Putting everything together, we get one-to-one correspondences
between the following four things, each considered up to the relevant sort of 
isomorphism:
\begin{itemize}
\item Dynkin diagrams
\item compact semisimple Lie algebras
\item connected and simply connected compact semisimple Lie groups
\item connected semisimple Lie groups with trivial center.
\end{itemize}
The really deep correspondence is between the first two items.    Then, for every
compact semisimple Lie algebra \(\g\), there are various connected Lie groups \(G\) having \(\g\) as their
Lie algebra.  They're all compact and semisimple, but some are covering spaces of others.
At the one extreme is the universal cover, which covers all the others: this is simply connected.
At the other extreme is the so-called \define{adjoint form}, which is covered \emph{by} all the
others: this has trivial center.     For example, \(\SU(N)\) is simply connected, but it has a nontrivial
center, consisting of all matrices that are an \(N\)th root of unity times the identity matrix.
If we mod out by this center we get the adjoint form, which is called the \define{projective special
unitary group} \(\mathrm{PSU}(N)\).

Furthermore, it turns out that the operation of taking products of compact semisimple
Lie groups corresponds to taking direct sums of their Lie algebras, which 
corresponds to taking disjoint unions of Dynkin diagrams.  So to get the 
``building blocks'' from which everything else can be built,
we only need to worry about the \emph{connected} Dynkin diagrams, which we
have completely classified:

\begin{center}
 \includegraphics[scale = 0.14]{dynkin_diagrams.png}
\end{center}

The compact semisimple Lie algebras coming from \emph{connected} Dynkin diagrams
are called \define{compact simple Lie algebras}.  But what are they actually like?     People have figured them out.   Amazingly, those that come in infinite series are all related to
rotations in real, complex or quaternionic vector spaces.   So let me give you a crash course
on those:

\begin{enumerate}
\item The real-linear transformations of \(\R^n\) that preserve its usual inner product
\[         \langle v , w \rangle = \sum_{i=1}^n v_i w_i \]
form a compact Lie group called the \define{orthogonal group} \(\O(n)\).  The subgroup consisting of transformations with determinant \(1\) is a compact Lie group called the \define{special orthogonal group} \(\SO(n)\).   Its Lie algebra is called  \(\mathfrak{so}(n)\), and it consists of \(n \times n\) real matrices that have trace zero and are minus their own transpose.  This is
a compact simple Lie algebra when \(n \ge 3\).  What is it when \(n = 1\) or \(2\)?
\item The complex-linear transformations of \(\C^n\) that preserve its usual inner product
\[         \langle v , w \rangle = \sum_{i=1}^n \overline{v}_i w_i \]
form a compact Lie group called the \define{unitary group} \(\U(n)\).  The subgroup consisting of transformations with determinant \(1\) is a compact Lie group called the \define{special unitary group} \(\SU(n)\).   Its Lie algebra is called  \(\mathfrak{su}(n)\), and it consists of \(n \times n\) complex matrices that have trace zero and are minus their own conjugate transpose.  This is
a compact simple Lie algebra when \(n \ge 2 \).    What is it when \(n = 1\)?
\item The quaternion-linear transformations of \(\H^n\) that preserve its usual inner product
\[         \langle v , w \rangle = \sum_{i=1}^n \overline{v}_i w_i \]
form a Lie group called the \define{quaternionic unitary group} \(\Sp(n)\).   (With Dieudonne's definition of the quaternionic determinant, all matrices in this group
have determinant \(1\).)    The Lie algebra of \(\Sp(n)\) is called \(\mathfrak{sp}(n)\), and
it consists of \(n \times n\) quaternionic matrices that have trace zero and are minus their own conjugate transpose.  This is a compact simple Lie algebra when \(n \ge 1\).
\end{enumerate}
The third item looks a lot like the first two, but you may be unfamiliar with the quaternions 
\(\H\).  For now I'll just say that there are three finite-dimensional associative algebras over 
\(\R\) equipped with a norm obeying 
\[                    |ab| = |a| |b| .\]
These are the real numbers \(\R\), the complex numbers \(\C\), and most excitingly the 
\define{quaternions} 
\[           \H = \{a + bi + cj + dk \, | \, a,b,c,d \in \R \} \] 
where multiplication is determined by the equations Hamilton carved into a wall on
on the 16th of October in 1843: 
\[      i^2 = j^2 = k^2 = i j k = -1 .\]
We can define conjugation for quaternions by
\[         q = a + bi + cj + dk \implies \overline{q} = a - bi - cj - dk \]
and then the norm is given by 
\[         |q| = \sqrt{qq^*} = \sqrt{q^* q} . \]  
The space \(\H^n\) acts like a ``quaternionic vector space'' even though \(\H\) is not a field: we
can multiply vectors on the left or the right by quaternions, and we say that \(T \maps \H^n \to \H^n\) is \define{quaternion-linear} if
\[      T(vq) = T(v) q  \qquad \textrm{ for all } v \in \H^n, q \in \H.\]
With this convention, left multiplication by any \(n \times n\) matrix of quaternions gives a quaternion-linear map \(T \maps \H^n \to \H^n\).   While the quaternions are a fascinating
subject, this is all you need to know for now.

Now we are ready to list the four infinite series of compact simple Lie algebras:

\begin{itemize}
\item \define{\(\mathrm{A}_n\)}: This Dynkin diagram gives the compact simple Lie algebra 
\(\mathfrak{su}(n+1)\).
\item
\define{\(\mathrm{B}_n\)}: This Dynkin diagram gives the compact simple Lie algebra 
\(\mathfrak{so}(2n+1)\).
\item
\define{\(\mathrm{C}_n\)}:  This Dynkin diagram gives the compact simple Lie algebra
 \(\mathfrak{sp}(n)\).
\item
\define{\(\mathrm{D}_n\)}:  This Dynkin diagram gives the compact simple Lie algebra 
\(\mathfrak{so}(2n)\).
\end{itemize}

These are called the  \define{classical} compact simple Lie algebras, and they would be
pretty easy to reinvent for yourself, or get interested in for all
sorts of reasons.    It may seem weird that 
\(\mathrm{SO}(2n)\) is so different from \(\mathrm{SO}(2n+1)\), but it's true!   For example, can put \(n\) orthogonal planes in \(\R^{2n}\), and by doing a rotation in each of these planes you get an element of \(\mathrm{SO}(2n)\) that fixes only the origin.   But in odd dimensions 
there's one dimension left over, so any rotation must fix some nonzero vector.

The remaining five compact simple Lie algebras are called \define{exceptional}, and
 they are much more
mysterious. They were only discovered when people like Killing and Cartan figured out 
the classification of simple Lie algebras. And as it turns out, they are all 
related to the octonions!   The octonions \(\Oct\) are the only 
finite-dimensional \emph{nonassociative} algebra over \(\R\) that is equipped with a norm obeying
\[            |ab| = |a| |b| .\]
They are an amazing freak of nature, which begets many other strange things.   I am somewhat
obsessed with them, but I won't say much about them here: see the references for more.

Here are the 5 exceptional compact simple Lie algebras.   It is quickest to describe
them using compact Lie groups.   I'll list them in order of dimension, not alphabetical order:

\begin{enumerate}
\item
\define{\(\mathrm{G}_2\)}:  This Dynkin diagram gives a 14-dimensional compact simple Lie algebra called \(\mathfrak{g}_2\).   The automorphism group of the octonions is a 
compact Lie group whose Lie algebra is \(\mathfrak{g}_2\).
\item
\define{\(\mathrm{F}_4\)}:  This Dynkin diagram gives a 52-dimensional compact simple Lie algebra called \(\mathfrak{f}_4\).  The octonions give a 16-dimensional projective plane
called \(\Oct\mathrm{P}^2\).  This is a Riemannian manifold, and its \define{isometry
group}---the group of diffeomorphisms preserving the Riemannian metric---is a compact
Lie group whose Lie algebra is \(\mathfrak{f}_4\).
\item
\define{\(\mathrm{E}_6\)}:  This Dynkin diagram gives a 78-dimensional compact simple Lie algebra called \(\mathfrak{f}_4\).  The octonions tensored with the complex numbers 
give a 32-dimensional projective plane called \((\C \otimes \Oct)\mathrm{P}^2\).  This is a Riemannian manifold, and its isometry group is a compact Lie group whose Lie algebra is 
\(\mathfrak{e}_6\).
\item
\define{\(\mathrm{E}_7\)}:  This Dynkin diagram gives a 133-dimensional compact simple Lie algebra called \(\mathfrak{e}_7\).  The octonions tensored with the quaternions 
give a 64-dimensional Riemannian manifold called \((\H \otimes \Oct)\mathrm{P}^2\),
even though it is not technically a projective plane.   Its isometry group is a compact
Lie group whose Lie algebra is \(\mathfrak{e}_7\).
\item
\define{\(\mathrm{E}_8\)}:  This Dynkin diagram gives a 248-dimensional compact simple Lie algebra called \(\mathfrak{e}_8\).  The octonions tensored with the octonions 
give a 128-dimensional Riemannian manifold called \((\Oct \otimes \Oct)\mathrm{P}^2\),
even though it is not technically a projective plane.   Its isometry group is a compact
Lie group whose Lie algebra is \(\mathfrak{e}_8\).
\end{enumerate}

A tragic fact is that currently nobody know how to get their hands on the 
Riemannian manifolds  \((\H \otimes \Oct)\mathrm{P}^2\) and 
\((\Oct \otimes \Oct)\mathrm{P}^2\) without building their isometry groups first, or
at the same time.   Thus the explanations we have given of \(\mathfrak{e}_7\) and \(\mathfrak{e}_8\), while true, are not as useful as we might like.   This is especially 
frustrating for \(\mathfrak{e}_8\), since this is the ``king'' of simple Lie algebras, connected
to many other amazing exceptional structures in mathematics.

We've listed the compact simple Lie algebras, but what about their Lie groups?  This is an important subject, since Lie algebras are ultimately just a tool for working with Lie groups.  But we need to be a bit careful, since a Lie algebra may be the Lie algebra of several nonisomorphic connected Lie groups.  Remember, taking a covering space of a Lie group does not change its Lie algebra.   For example, the Dynkin diagram coincidence \(\mathrm{A}_1 \cong \mathrm{B}_1\) implies that \(\mathfrak{su}(2) \cong \mathfrak{so}(3)\), but the Lie group \(\SU(2)\) is not isomorphic to \(\SO(3)\): it's a double cover of \(\SO(3)\).   

People usually take a relaxed attitude and call \emph{any} connected Lie group whose Lie algebra is a compact simple Lie algebra a \define{compact simple Lie group}.  This is true even though such a group may have nontrivial normal subgroups, so it is not simple in the usual sense of group theory.   For example \(\SU(N)\) has an \(n\)-element normal subgroup, its center, containing the matrices that equal an \(N\)th root of unity times the identity matrix.   Yet we call it a compact simple Lie group.    If you don't like this, fear not: if you take any compact simple Lie group and mod out by its center, you get another compact simple Lie group with the same Lie algebra, which doesn't have any nontrivial normal subgroups.  This is the so-called ``adjoint form'', which I mentioned earlier.

With this terminology in place, we get a one-to-one correspondence
between these four things, each considered up to the relevant sort of 
isomorphism:
\begin{itemize}
\item connected Dynkin diagrams
\item compact simple Lie algebras
\item simply connected compact simple Lie groups
\item compact simple Lie groups without nontrivial normal subgroups.
\end{itemize}
Among the the classical compact simple Lie groups it turns out that \(\SU(n)\) and \(\Sp(n)\) are simply connected, while \(\SO(n)\) is not: for \(n \ge 3\), which is all that matters here, \(\SO(n)\) has a simply connected double cover called the \define{spin group} \(\mathrm{Spin}(n)\).  This group is very important in physics and also differential geometry, since it has important representations called ``spinors'' that are not representations of \(\SO(n)\).   I have written
quite a bit about the exceptional compact simple Lie groups, but here I'll just mention that any connected Lie group with Lie algebra \(\mathfrak{g}_2, \mathfrak{f}_4, \mathfrak{e}_6, \mathfrak{e}_7\) or \(\mathfrak{e}_8\) is called \(\mathrm{G}_2, \mathrm{F}_4, \mathrm{E}_6, \mathrm{E}_7\) or \(\mathrm{E}_8\).  This is uncreative and a bit ambiguous, but that's that way it is.

\subsubsection*{Simply-laced Dynkin diagrams}

We've seen that any Dynkin diagram with \(n\) dots describes a collection of
vectors \(v_1, \dots, v_n \in \R^n\) such that:
\begin{enumerate}
\item
integer linear combinations of these vectors form a lattice \(L \subset \R^n\)
\item 
reflections through these vectors generate a finite reflection group \(\Gamma\)
\item
the action of \(\Gamma\) on \(\R^n\) preserves the lattice \(L\).
\end{enumerate}
Unfortunately, we can't always choose all the vectors \(v_i\) to have the same
length.   But we've seen that this problem doesn't happen when all these vectors 
are at an angle of \(\pi/2\) or \(\pi/3\) from each other.   This happens when our
Dynkin diagram is \define{simply laced}: all its edges are unlabeled.

Here are all the connected simply-laced Dynkin diagrams:
\begin{itemize}
\item
\(\mathrm{A}_n\), which has \(n\) dots like this:
 \[
  \begin{tikzpicture}
    \draw[very thick] (0,0) node{$\sbullet$} to (1,0) node{$\sbullet$} to (2,0) node{$\sbullet$} to (3,0) node{$\sbullet$};
  \end{tikzpicture}
\]
\item
 \(\mathrm{D}_n\), which has \(n\) dots, where we can assume \(n > 3\) since \(\mathrm{D}_n \cong \mathrm{A}_n\) for \(n = 1,2,3\): 
 \[
  \begin{tikzpicture}
    \draw[very thick] (0,0) node{$\sbullet$} to (1,0) node{$\sbullet$} to (2,0) node{$\sbullet$} to (3,0) node{$\sbullet$};
    \draw[very thick] (3,0) to (4,1) node{$\sbullet$};
    \draw[very thick] (3,0) to (4,-1) node{$\sbullet$};
  \end{tikzpicture}
\]  
\item \(\mathrm{E}_6\), \(\mathrm{E}_7\), and \(\mathrm{E}_8\): \[
  \begin{gathered}
    \begin{tikzpicture}
      \draw[very thick] (0,0) node{$\sbullet$} to (1,0) node{$\sbullet$} to (2,0) node{$\sbullet$} to (3,0) node{$\sbullet$} to (4,0) node{$\sbullet$};
      \draw[very thick] (2,0) to (2,-1) node{$\sbullet$};
    \end{tikzpicture}
\\\begin{tikzpicture}
      \draw[very thick] (0,0) node{$\sbullet$} to (1,0) node{$\sbullet$} to (2,0) node{$\sbullet$} to (3,0) node{$\sbullet$} to (4,0) node{$\sbullet$} to (5,0) node{$\sbullet$};
      \draw[very thick] (3,0) to (3,-1) node{$\sbullet$};
    \end{tikzpicture}
  \\\begin{tikzpicture}
      \draw[very thick] (0,0) node{$\sbullet$} to (1,0) node{$\sbullet$} to (2,0) node{$\sbullet$} to (3,0) node{$\sbullet$} to (4,0) node{$\sbullet$} to (5,0) node{$\sbullet$} to (6,0) node{$\sbullet$};
      \draw[very thick] (4,0) to (4,-1) node{$\sbullet$};
    \end{tikzpicture}
  \end{gathered}
\]   
\end{itemize}
These diagrams are ubiquitous in mathematics.   But before getting into that
I should describe the corresponding lattices more explicitly, to
make it clear how simple they really are.

So, what are the \(\mathrm{A}, \mathrm{D}\), and \(\mathrm{E}\) lattices?

\begin{itemize}
\item
\define{\(\mathrm{A}_n\)}:  We can describe the \(\mathrm{A}_n\) lattice as 
the set of all \((n+1)\)-tuples of integers \((x_1,\ldots,x_{n+1})\) such that
\[x_1+\cdots+x_{n+1}=0.\] 
It's a fun exercise to show that \(A_2\) is a
\(2\)-dimensional hexagonal lattice, the sort of lattice you use to pack
pennies as densely as possible. Similarly, \(\mathrm{A}_3\) gives a standard way
of packing cannon balls, which is the densest lattice packing of equal-sized
spheres in 3 dimensions.  A much harder fact, due to Hales, is that no non-lattice
packing of equal-sized spheres can beat the density of the \(\mathrm{A}_3\) lattice.
\item
\define{\(\mathrm{D}_n\)}: We can describe the \(\mathrm{D}_n\) lattice
as the set of all \(n\)-tuples of integers \((x_1,\ldots,x_n)\) such that
\[x_1+\cdots+x_n\quad\text{is even}.\] Or, if you like, you can imagine
taking an \(n\)-dimensional checkerboard, coloring the cubes alternately
red and black, and taking the center of each red cube. In four
dimensions, \(D_4\) gives a denser packing of spheres than \(A_4\); in
fact, it gives the densest lattice packing possible. Moreover, \(D_5\)
gives the densest lattice packing of in dimension 5. However, in
dimensions 6, 7, and 8, the \(\mathrm{E}_n\) lattices give the densest lattice
packings.  In fact Viasovska showed that in 8 dimensions, no non-lattice
packing of equal sized spheres can beat the density of the \(\mathrm{E}_8\)
lattice!
\item
\define{\(\mathrm{E}_6, \mathrm{E}_7, \mathrm{E}_8 \)}:  We can describe 
the \(\mathrm{E}_8\) lattice as the set of 8-tuples
\((x_1,\ldots,x_8)\) such that the \(x_i\) are either all integers or all
integers plus \(1/2\) and
\[x_1+\cdots+x_8\quad\text{is even}.\] 
Each point in this lattice has 240 nearest neighbors.  For example, the
nearest neighbors of the origin have length \(\sqrt{2}\), and you can check
there are 240 of them.

Alternatively, if you take the \(\mathrm{D}_8\) lattice and use it to pack equal-sized spheres 
that just touch each other, there is actually just enough room to slip in another 
\(\mathrm{D}_8\)  lattice of equal-sized spheres in the remaining space, doubling the density!
And if you do this, your spheres will be centered at points in the \(\mathrm{E}_8\)
lattice.

Once you have \(\mathrm{E}_8\) in hand, you can get its little pals
\(\mathrm{E}_7\) and \(\mathrm{E}_6\) as follows. To get
\(\mathrm{E}_7\), just take all the vectors in \(\mathrm{E}_8\) that are
perpendicular to one lattice vector of length \(\sqrt{2}\).   
To get \(\mathrm{E}_6\), find a copy of the lattice \(A_2\) in \(\mathrm{E}_8\)
generated by 2 vectors of length \(\sqrt{2}\),  and then take all the vectors in 
\(\mathrm{E}_8\) perpendicular to everything in that copy of \(A_2\).
\end{itemize}

The A, D and E Dynkin diagrams show up in many places throughout 
mathematics, in a spooky sort of way.   Let me sketch three of the most famous.

First, Witt's theorem says that the A, D, and E lattices and their
direct sums are the only integral lattices having a basis consisting of
vectors \(v\) with \(\|v\|^2 = 2\). Here a lattice is \define{integral} if
the dot product of any two vectors in it is an integer. In fact, any
integral lattice having a basis consisting of vectors with \(\|v\|^2\)
equal to \(1\) or \(2\) is a direct sum of copies of A, D, and E
lattices and the integers, thought of as a \(1\)-dimensional lattice.

Second, a \define{quiver} is just some dots with arrows between them.   A 
\define{representation} of a quiver is a way of assigning a finite-dimensional
complex vector space to each dot and a linear map between these vector spaces
to each arrow.  There's an obvious category of representations \(\mathrm{Rep}(Q)\)
of any quiver \(Q\).   Gabriel proved an astounding result about these categories 
\(\mathrm{Rep}(Q)\). We say a quiver \(Q\) has \define{finite representation type}
if \(\mathrm{Rep}(Q)\) has finitely many isomorphism classes of 
\define{indecomposable} objects: objects that aren't direct sums of others.  And, it turns 
out the quivers of finite representation type are just those coming from simply-laced Dynkin diagrams!

Actually, for this to make sense, you need to take your Dynkin diagram and turn it into a 
quiver by putting arrows along the edges. If you have a simply-laced Dynkin diagram, you get 
a quiver of finite representation type no matter which way you let the arrows point. 

Third, there is a cool relationship between the ADE
diagrams and the symmetry groups of the Platonic solids, called the
\define{McKay correspondence}.  Here is one way to get it.  First, take the
rotational symmetry group of a Platonic solid, not including reflections,
or more generally any finite subgroup \(G\) of \(\mathrm{SO}(3)\). 
Since \(\mathrm{SO}(3)\) has
\(\mathrm{SU}(2)\) as a double cover, you can get a double cover of
\(G\), say \(\widetilde{G}\), sitting inside \(\mathrm{SU}(2)\). 
Since \(\widetilde{G}\) is finite, it has finitely many irreducible 
representations on complex vector spaces (up to isomorphism).  
Draw a dot for each of these.  One comes from the obvious
representation of \(\SU(2)\) on \(\C^2\).   When you tensor this one 
with any other irreducible representation \(R\), you get a direct sum of
irreducible representations.     Draw one line from the dot for \(R\) to another dot
for each time that other irreducible representation appears in your direct sum.
What do you get? 

You get an ``affine Dynkin diagram'', which is like a usual Dynkin diagram 
but with an extra dot thrown in---corresponding to the trivial rep of \(\widetilde{G}\).   
And if you throw out that extra dot, you get a simply laced Dynkin diagram!
In fact you get all all the connected simply laced Dynkin diagrams this way!  

The correspondence goes like this:
\begin{itemize}
\item
\define{\(\mathrm{A}_n\)}: this corresponds to the cyclic group
sitting inside \(\SO(n)\) as the rotational symmetries of a regular \(n\)-gon, where
we don't let ourselves flip this polygon over.
\item 
\define{\(\mathrm{D}_n\)}: this corresponds to the dihedral group sitting inside 
 \(\SO(n)\) as the symmetries of a regular \(n\)-gon, where we \emph{do} let
 ourselves flip this polygon over.
 \item 
\define{\(\mathrm{E}_6\)}: this corresponds to the rotational symmetry group of the regular
tetrahedron.
\begin{center}
\includegraphics[scale=0.2]{tetrahedron.jpg}  
\end{center}
\item
\define{\(\mathrm{E}_7\)}: this corresponds to the rotational symmetry group of the regular
cube, or octahedron.
\begin{center}
\includegraphics[scale=0.2]{hexahedron.jpg}  \qquad \qquad
\includegraphics[scale=0.2]{octahedron.jpg} 
\end{center}
\item
\define{\(\mathrm{E}_8\)}: this corresponds to the rotational symmetry group of the regular dodecahedron, or icosahedron.
\begin{center}
\includegraphics[scale=0.2]{dodecahedron.jpg}  \qquad \qquad
\includegraphics[scale=0.2]{icosahedron.jpg} 
\end{center}
\end{itemize}
We saw another relation between Platonic solids and Coxeter diagrams near the start
of this paper, but that one made sense.  This one is black magic.    

\subsubsection*{References}

I have skimmed a lot of material but not explained it in detail. The lectures I gave at the
University of Edinburgh, based on these notes, may help a bit:
\begin{itemize}
\item 
John C. Baez, \href{https://math.ucr.edu/home/baez/twf/}{Talks on This Week’s Finds in 
Mathematical Physics}, Lectures 4--8.
\end{itemize}
But there are still many important details missing!   Here are some ways to learn more.

For a gentle introduction to symmetry groups in 2 and 3 dimensions and a valuable
warmup for higher-dimensional considerations, see:
\begin{itemize}
\item
D.\ L.\ Johnson, \textsl{Symmetries}, Springer, Berlin, 2001.
\end{itemize}
For the classification of finite reflection groups in terms of Coxeter diagrams, see:
\begin{itemize}
\item
 James E.\ Humphreys, \textsl{Finite Reflection Groups and Coxeter Groups},
Cambridge U.\ Press, Cambridge, 1990.
\end{itemize}
For their appearance as symmetry groups of polytopes, try this endlessly entertaining book:
\begin{itemize}
\item 
Harold Scott Macdonald Coxeter, \textsl{Regular Polytopes}, Dover, Mineola, 1974.
\end{itemize}
For a huge amount of information on lattices, see this book:
\begin{itemize}
\item
 John Horton Conway and Neil James Alexander Sloane,
  \textsl{Sphere Packings, Lattices and Groups}, Springer, Berlin, 2013.
\end{itemize}
Gian-Carlo Rota said of this book, ``This is the best survey of the best work in the best fields of combinatorics written by the best people.  It will make the best reading by the best students interested in the best mathematics that is now going on.''
Chapter 4 discusses the lattices coming from Dynkin diagrams.

For the classification of compact Lie groups, see:
\begin{itemize}
\item
  John Frank Adams, \textsl{Lectures on Lie Groups}, Benjamin, New York,
  1969.
\end{itemize}
This is nice because it uses topology and handles compact Lie groups that are not semisimple on a more or less equal footing with the semisimple ones.   There are many texts that take a more
algebraic approach to classifying semisimple Lie algebras and/or Lie groups:
Here are a few:
\begin{itemize}
\item  Daniel Bump, \textsl{Lie Groups}, Springer, Berlin, 2004. 
\item William Fulton and Joe Harris, \textsl{Representation Theory --- a First Course}, Springer, Berlin, 1991. 
\item Sigurdur Helgason, \textsl{Differential Geometry, Lie Groups, and 
Symmetric Spaces}, Academic Press, New York, 1979.
\item James Humphreys, \textsl{Introduction to Lie Algebras and Representation
Theory}, Springer, Berlin, 2012.
\end{itemize}
This is a useful reference:
\begin{itemize}
\item Nicolas Bourbaki, \emph{Lie Groups and Lie Algebras, Chapters 7--9}, Springer, Berlin, 1975.
\end{itemize}
Often such books start by
classifying semisimple Lie algebras over \(\C\), and then show each one is
the complexification of a unique compact semisimple Lie algebra.  
Each such complex Lie
algebra is also the Lie algebra of a unique connected and simply connected 
\define{complex Lie group}, meaning a group in the
category of complex manifolds and holomorphic maps.   So, there is actually a one-to-one correspondence
between six things, each considered up to isomorphism:
\begin{itemize}
\item connected Dynkin diagrams
\item compact simple Lie algebras
\item connected and simply connected compact simple Lie groups
\item connected compact simple Lie groups with trivial center (or equivalently, without nontrivial normal subgroups)
\item complex simple Lie algebras
\item connected and simply connected complex simple Lie groups.
\end{itemize}
In the body of this paper I assume connectedness as part of the definition of 
``simple Lie group'', but now I'm making it explicit.  

For more on the exceptional Lie groups and their connection to octonions, see:
\begin{itemize}
\item John Frank Adams, \textsl{Lectures on Exceptional Lie Groups}, eds.\
Zafer Mahmud and Mamoru Mimura, U.\ Chicago Press, Chicago, 1996.
\item John C.\ Baez, \href{https://arxiv.org/abs/math/0105155}{The octonions}, 
\textsl{Bull.\ Amer.\ Math.\ Soc.\ } \textbf{39} 
(2002), 145--205.  Errata in \textsl{Bull.\ Amer.\ Math.\ Soc.\ } \textbf{42} (2005), 213. 
\item Ichiro Yokota, \href{https://arxiv.org/abs/0902.0431}{Exceptional Lie
groups}.
\end{itemize}

For more on the ubiquitous appearance of Coxeter and Dynkin diagrams, and especially
the ADE Dynkin diagrams, see:
\begin{itemize}
\item
  M. Hazewinkel, W. Hesselink, D. Siermsa, and F. D. Veldkamp,
 \href{http://math.ucr.edu/home/baez/hazewinkel\_et\_al.pdf}{The
  ubiquity of Coxeter--Dynkin diagrams (an introduction to the ADE
  problem)}, \textsl{Niew.\ Arch.\ Wisk.} \textbf{25} (1977), 257--307.
\item 
  John McKay,  \href{http://math.ucr.edu/home/baez/ADE.html}{A rapid introduction to ADE theory}.
\item
  Vladimir I.\ Arnol'd, ``Problems of Present Day Mathematics'' in \emph{Mathematical
  Developments Arising from Hilbert's Problems}, ed.~F. E. Browder,
  Proc. Symp. Pure Math. \textbf{28}, AMS,
  Providence, Rhode Island, 1976.
\end{itemize}
Arnol'd lists a lot of important math problems, following up on Hilbert's
famous turn-of-the-century listing of problems. Problem VIII in this
book is the ``ubiquity of ADE classifications''. 

Witt showed that all the integral lattices generated by vectors \(v\) with 
\(v \cdot v = 2\) come from simply laced Dynkin diagrams in this paper:
\begin{itemize}
\item Ernst Witt, Spiegelungsgruppen und Aufz\"ahlung halbeinfacher Liescher Ringe, 
\textsl{Abhandlungen aus dem Mathematischen Seminar der Universit{\"a}t Hamburg} \textbf{14} (1941), 289--322.
\end{itemize}
There should be some more useful modern reference!
For a proof of Gabriel's theorem that only simply laced Dynkin diagrams give quivers of 
tame representation type, see:
\begin{itemize}
\item
 William Crawley-Boevey, \href{https://www.math.uni-bielefeld.de/~wcrawley/quivlecs.pdf}{Lectures on representations of quivers}.
\end{itemize}
For more discussions of Gabriel's theorem, see:
\begin{itemize}
\item 
Harm Derksen and Jerzy Weyman, \href{http://www.ams.org/notices/200502/fea-weyman.pdf}{Quiver representations}, \textsl{Notices Amer.\ Math.\ Soc.\ }\textbf{52} (2005), 200--206. 
\item 
Idun Reiten, Dynkin diagrams and the representation theory of algebras,
\textsl{Notices Amer.\ Math.\ Soc.} \textbf{44} (1977), 546--558.
\item
Alistair Savage, \href{http://arxiv.org/abs/arXiv:math/0505082}{Finite-dimensional algebras and quivers}.
\end{itemize}
For the McKay correspondence, see:
\begin{itemize}
\item
  John McKay, Graphs, singularities and finite groups, in
  \textsl{Proc.\ Symp.\ Pure Math.} vol \textbf{37}, AMS, Providence, 
  Rhode Island,  1980, pp.\ 183--186.
\item
  D.\ Ford and John McKay, Representations and Coxeter graphs, in \textsl{The
  Geometric Vein} Coxeter Festschrift (1982), Springer, Berlin,
  pp.\ 549--554.
\item
  Pavel Etinghof and Mikhail Khovanov,   \href{https://arxiv.org/abs/hep-th/9408078}{Representations of tensor categories and Dynkin diagrams}.
 \item
 Joris van Hoboken, \textsl{\href{
http://math.ucr.edu/home/baez/joris\_van\_hoboken\_platonic.pdf}{Platonic solids, binary polyhedral groups, Kleinian singularities and Lie algebras of type \(A\),\(D\),\(E\)}},
  Master's Thesis, University of Amsterdam, 2002.
\item
Klaus Lamotke, \textsl{Regular Solids and Isolated Singularities}, Vieweg \& Sohn, Braunschweig, 1986.
\item Peter Slodowy, Platonic solids, Kleinian singularities, and Lie groups, in 
\textsl{Algebraic Geometry}, Springer, Berlin, 1983, pp.\ 102--138.
\end{itemize}
As you can see, some of these describe the connection between the ADE diagrams
and singularity theory, which I did not touch on here.

\vskip 1em
The green images of Platonic solids were created by \href{https://commons.wikimedia.org/wiki/User:Kjell_Andre}{Kjell Andr\'e} and converted into SVG files by \href{https://commons.wikimedia.org/wiki/User:DTR}{DTR}.  They are available at Wikicommons
under a CC BY-SA 3.0 license.  The images of tilings and the fancy images of Dynkin 
diagrams were created by \href{https://commons.wikimedia.org/wiki/User:Nonenmac}{R.\ A.\ Nonemacher} and are available at Wikicommons under a CC BY-SA 4.0 license.   The other images were either created by myself, or I've been unable to track down their source.

\end{document}